\documentclass[reqno]{amsart}
\usepackage{amsthm, amssymb, amsfonts}
\usepackage{lmodern}
\usepackage{color}
\usepackage{hyperref}
\usepackage[T5]{fontenc}
\usepackage{amscd,amssymb}
\usepackage[v2,cmtip]{xy}
\usepackage{mathrsfs}
\usepackage[top=3.0cm, bottom=2.5cm, left=3.7cm, right=3.7cm] {geometry}
\theoremstyle{plain}

\theoremstyle{definition}

\theoremstyle{plain}

\theoremstyle{definition}

\allowdisplaybreaks

\begin{document} 
\title[On an article published in the racsam]
{On an article published in the racsam}

\author{Nguy\~\ecircumflex n Sum}
\address{Department of Mathematics and Applications, S\`ai G\`on University, 273 An D\uhorn \ohorn ng V\uhorn \ohorn ng, District 5, H\`\ocircumflex\ Ch\'i Minh city, Viet Nam}
 
\email{nguyensum@sgu.edu.vn}

\footnotetext[3]{The main contents of this Note have been submitted to Mathematical Reviews of AMS and to the managing Editors of the RACSAM.}


\begin{abstract}
In this Note, we provide the comments on the paper \cite{phc}  [\textit{A note on the hit problem for the polynomial algebra in the case of odd primes and its application}] which was published in the RACSAM [Rev. Real. Acad. Cienc. Exactas. Fis. Nat. Ser. A-Mat. 118, 22 (2024)]. We will show that some of the results in this paper are not new and others are false.
\end{abstract}
\maketitle

\section{Introduction}\label{s1} 

Let $p$ be an odd prime number and $\mathscr A_p$ be the mod-$p$ Steenrod algebra. Denote by $P_h = \mathbb F_p[t_1,t_2,\ldots ,t_h]$ the polynomial algebra over the field $\mathbb F_p$ with $p$ elements in $h$ generators $t_1,t_2,\ldots ,t_h$, each of degree 2. This algebra is regarded as an $\mathscr A_p$ module. Let $GL(h,\mathbb F_p)$ be the general linear group of rank $h$ over $\mathbb F_p$ which acts regularly on $P_h$. 

In the paper \cite{phc}, the author presented some results on $\mathscr A_p$-generators of $P_h$ and applied them to study the third mod-$p$ Singer algebraic transfer which was defined in Crossley \cite{cros2}. The main results of the paper are theorems numbered 2.1 and from 2.4 to 2.8. However, some of the results in the paper \cite{phc} are not new and others are false. We use the notations as in \cite{phc}.

\medskip
1. It is well known that $P_1$ is an $\mathscr A_p$ submodule of $P_2$. A monomial basis for $\mathbb F_p\otimes_{\mathscr A_p}P_1$ is a subset of the one for $\mathbb F_p\otimes_{\mathscr A_p}P_2$, hence Theorem 2.1 is only a special case of Theorem 1.1 on Page 172 of the paper Crossley \cite{cros}. It is not new.

\medskip
2. On Page 5, the author presented an example which substantially overlaps with the one in Crossley \cite[Page 175]{cros}. However, the author did not cite the paper \cite{cros} of M. D. Crossley.

\medskip
3. There are serious mistakes in the proof of Theorem 2.4, so I think that this theorem is false.

\medskip
4. The part ``$g^p$ is hit'' of Theorem 2.5 is not new. It is only the unstability of $\mathscr A_p$ module $P_h$ which is listed in Tannay-Oner \cite[Prop. 1.4 on Page 16]{tan}.

\medskip
5. The statement of Theorem 2.6 is not correct. Moreover, the result of this theorem is trivial.

\medskip
6. Theorem 2.7 is stated without proof. However, its statement is also not correct. I think that the author did not prove this theorem. Hence, its accuracy is not determined.

\medskip
7. The proof of Theorem 2.8 is refused with the following reasons: 

- The computation of $[(\mathbb F_p\otimes _{\mathscr A_p}H^*(V;\mathbb F_p))_{n^{(s)}}]^{GL(3,\mathbb F_p)}$ is based on a basis of the vector space $(\mathbb F_p\otimes _{\mathscr A_p}H^*(V;\mathbb F_p))_{n^{(s)}}$. However, one doesn't know any basis for this space yet.

- This proof is based on Theorems 2.4, 2.6, 2.7 but the accuracy of these results are not determined.

- There are infinite many odd prime numbers but computer programs based on algorithm of MAGMA can only check the space $(\mathbb F_p\otimes _{\mathscr A_p}H^*(V;\mathbb F_p))_{n^{(s)}}$ for a finite number of prime $p$. 

Thus, the results of Theorem 2.8 and Corollary 2.11 have not yet been accurately determined. 

\medskip
Below is the detailed comments for the above items.
 
\section{The detailed comments}\label{s2}
\setcounter{equation}{0}

\begin{enumerate}	
	\item[1.]  In \cite[Page 172]{cros}, Crossley stated Theorem~1.1 as follows:
	
	\medskip\noindent
	\textbf{Theorem 1.1.} \textit{As far as it goes, Table 1 gives a monomial basis for $M^*(2)$.}
	
	\medskip
	\centerline{Table 1}
	
	\medskip
	\rightline{\begin{tabular}{ll}
			\hline
			Degree, $n$ & Basis for $M^n(2)$\\
			\hline
			$\leq 2(p-2)$, $n$ even &$\{\textcolor{red}{x^iy^{n/2-i}}\,|\, 0\leq i \leq n/2\}$\\
			\hline
			$2(((i + 1)p + j + 1)p^s - 2)$&$\{x^{(k+1)p^s-1}y^{n/2-(k+1)p^{s}+1}\, |\, \min(i + 1, j) \leq k \leq p - 1\}$\\
			$0 \leq i,\, j \leq p - 1,\, s \geq 0$&$\cup \{\textcolor{red}{x^{(k+1)p^{s+1}-1}y^{n/2-(k+1)p^{s+1}+1}}\, |\, 0 \leq k \leq i\}$\\
			\hline
			$2(((i + 1)p^r + j + 1)p^s - 2)$&$\{x^{(k+1)p^{s+1}-1}y^{n/2-(k+1)p^{s+1}+1}\, |\, 1 \leq k \leq p - 1\}$\\
			$1 \leq i,\, j + 1 \leq p - 1,$&$\cup\{x^{(j+1)p^s-1}y^{(i+1)p^{r+s}-1},\, \textcolor{red}{x^{(i+1)p^{r+s}-1}y^{(j+1)p^s-1}}\}$\\
			$r \geq 2,\, s \ge 0$&\\
			\hline
			$2((p^2 + ip + j + 1)p^s - 2)$& $\{x^{(k+1)p^{s+1}-1}y^{n/2-(k+1)p^{s+1}+1}\, |\, i \leq k \leq j\}$\\
			$1 \leq i \leq j \leq p - 2,\, s \geq 0$&\\
			\hline
		\end{tabular}\ \ } 
	
	\bigskip\noindent
	\textbf{Statement 1.2.} \textit{In the degrees, $n$, not dealt by Table 1,  $M^*(2)$ is $0$.}
	
	\medskip
	Note that $x=t_1,\, y = t_2$ and $M^*(2)=\mathbb F_p\otimes_{\mathscr A_p}P_2$. A monomial basis of $M^*(1) = \mathbb F_p \otimes_{\mathscr A_p} P_1$ consists of monomials in the basis for $M^*(2)$ whose power of $y=t_2$ is 0. Hence, from Theorem 1.1 one gets (the red monomials with power 0 for $y$).
	
	\medskip\noindent
	\textbf{Corollary.} \textit{Table 2 gives a monomial basis for $\mathbb F_p \otimes_{\mathscr A_p} P_1 = M^*(1)$.}
	
	\medskip
	\centerline{Table 2}
	
	\medskip
	\centerline{\begin{tabular}{ll|l}
			\hline
			Degree, $n$ & Basis for $M^n(1)$&Note\\
			\hline
			$\leq 2(p-2)$, $n$ even &$\{x^{n/2}\}$&The first red monomial with $i =  n/2$\\
			\hline
			$2(i + 1)p-2 $&$\{x^{(i+1)p-1}\}$&The second red monomial with \\
			$0\leq i \leq p-1$& &$s=j=0, k=i$\\
			\hline
			$2(i + 1)p^r-2$& $\{x^{(i+1)p^{r}-1}\}$&The last red monomial with\\
			$1 \leq i \leq p-1,\, r \geq 2$& &$s = j = 0$\\
			\hline
	\end{tabular}} 
	
	\medskip
	Thus, Theorem 2.1 in this paper overlaps with this corollary and Statement~1.2. Hence, it is not new.
	
	\medskip
	\item[2.] In \cite[Page 175]{cros}, Crossley presented an example as follows:
	
	\medskip
	``For example, let $p = 5$ and consider the element $x^4y^5 + x^8y : \mathscr P^1(x^4y^5 + x^8y) = (4x^8y^5 + 5x^4y^9) + (8x^{12}y + x^8y^5)=3x^{12}y$. Thus $x^{12}y$ is hit by $\mathscr P^1$.''
	
	\medskip
	On Page 5 of the paper \cite{phc}, the author presented an example as follows:
	
	\medskip
	``For instance, with $p = 5$ and the monomial
	$f = 3t_1^{12}t_2 \in P_2 = \mathbb F_5[t_1,t_2]$, we see that $f$ is in the image of $\mathscr A_5$ since $f = \mathscr P^1(t_1^4t_2^5+t_1^8t_2)$.''
	
	\medskip
	This example substantially overlaps with the one of Crossley. However, the author did not cite the paper \cite{cros} of M. D. Crossley.
	
	\medskip
	\item[3.] The proof of Theorem 2.4 is presented on Page 5 as follows:
	
	``We use induction on $h$. Obviously, the theorem is true for
	$h=1$. Suppose that the theorem is true for
	$P_{h-1}$. We write $f= t f_1 +f_2$, in which $t=t_h,\, f_1 \in P_h,\, \textcolor{red}{\deg(f_1)=n-2(p-1)}$ and $f_2 \in P_{h-1}$, $\deg(f_2) = n$. Then by the Cartan formula, we get
	$$0 = \mathscr P^1(f) = \mathscr P^1(tf_1) + \mathscr P^1(f_2) = t\mathscr P^1(f_1) + t^pf_1 + \mathscr P^1(f_2).$$
	
	As it is known, every term of $\mathscr P^k(f_2)\, (k \geqslant 0)$  involves exactly the same variables as $f_2$ does. So, since
	$f_2$ is independent of $t$, setting $t=0$, we get $\mathscr P^1(f_2) =0$ and so $t\mathscr P^1(f_1) + t^pf_1 = 0$. These, together with the inductive hypothesis, yield that $\textcolor{red}{\mathscr P^1(f_1) = tf_1}$ and $f_2 = \mathscr P^1(g_2)$ for some $g_2 \in P_{h-1}$, $\deg(g_2) = n - 2(p-1)$. Thus, $\textcolor{red}{f= t f_1 +f_2 = \mathscr P^1(f_1) + \mathscr P^1(g_2) = \mathscr P^1(f_1+g_2)}$. This completes the proof.'' 
	
	\medskip
	This proof is false because the author wrote $f = tf_1+f_2$ with $t = t_h$, hence the relation $\deg f_1 = n - 2(p-1)$ is false because $\deg t =2$, $n = \deg f$ and $\deg f_1 = n-2>n - 2(p-1)$. Since $p$ is odd, the relation $t\mathscr P^1(f_1)+ t^pf_1 =0$ implies $\mathscr P^1(f_1) = - t^{p-1}f_1 \ne tf_1$ for all $f_1 \ne 0$. Hence, the relations  $\mathscr P^1(f_1) = tf_1$ and
	$f= tf_1 + f_2 =\mathscr P^1(f_1 + g_2)$ are false. 
	
	\medskip
	\item[4.] Theorem 2.5 is only an easy consequence of the unstability of $\mathscr A_p$ module $P_h$ and the Cartan formula, however there are many errors in its proof (some of them are corrected in the Correction of \cite{phc}). Moreover, for any homogeneous element $g\in P_h$, one has $g^p = \mathscr P^{\frac 12\deg g}(g)$ (see Proposition~1.4 in the paper Tannay-Oner \cite[Page 16]{tan}). Hence, the condition $\deg g = 2p^j$ is unnecessary.
	
	\medskip
	\item[5.] The author stated Theorem 2.6 on Page 6 as follows:
	
	\medskip\noindent
	\textbf{Theorem 2.6.} \textit{For positive integers $t$, $q_i,\, 1\leqslant i \leqslant h$, let $\varphi:P_h \to P_h$ be the linear map defined by $\varphi (f) = t_1^{q_1p^{t+1}}t_2^{q_2p^{t+1}}\ldots t_h^{q_hp^{t+1}}f^p$ for all $f\in P_h$. Then if $f$ is in the image of $\mathscr A_p$ via $f = \sum_{0\leq k \leq t-1}\mathscr P^k(f_k)$ for some $f_k \in P_h$, then then $\varphi$ induces a homomorphism
		$$(\mathbb F_p\otimes_{\mathscr A_p} P_h)_n \to (\mathbb F_p\otimes_{\mathscr A_p} P_h)_{pn+2p^{t+1}\sum_{1\leqslant i \leqslant h}q_i} $$ for all $n \geq 0$.}
	
	The statement of this theorem is not correct because the map $\varphi$ induces a homomorphism if and only if $\varphi(f)$ is hit for all hit polynomials $f$. Hence, this theorem must be restated as follows:
	
	\medskip
	\textit{For positive integers $t$, $q_i,\, 1\leqslant i \leqslant h$, and let $\varphi:P_h \to P_h$ be the linear map defined by $\varphi (f) = t_1^{q_1p^{t+1}}t_2^{q_2p^{t+1}}\ldots t_h^{q_hp^{t+1}}f^p$ for all $f\in P_h$. If $f$ is hit, then so is $\varphi(f)$. Consequently, $\varphi$ induces a homomorphism
		$$(\mathbb F_p\otimes_{\mathscr A_p} P_h)_n \to (\mathbb F_p\otimes_{\mathscr A_p} P_h)_{pn+2p^{t+1}\sum_{1\leqslant i \leqslant h}q_i} \mbox{ for all } n \geq 0.$$}
	However, by Theorem 2.5,
	$$\varphi (f) = t_1^{q_1p^{t+1}}t_2^{q_2p^{t+1}}\ldots t_h^{q_hp^{t+1}}f^p = (t_1^{q_1p^{t}}t_2^{q_2p^{t}}\ldots t_h^{q_hp^{t}}f)^p$$
	is hit for any polynomial $f\in P_h$. So, $\varphi$ induces a trivial homomorphism from the space $(\mathbb F_p\otimes_{\mathscr A_p} P_h)_n$ to the space $(\mathbb F_p\otimes_{\mathscr A_p} P_h)_{pn+2p^{t+1}\sum_{1\leqslant i \leqslant h}q_i}$. Therefore, this theorem is very trivial.
	
	\medskip
	\item[6.] Theorem 2.7 is stated without proof:
	
	\medskip\noindent
	\textbf{Theorem 2.7.} \textit{For positive integers $t \geq 2$, $q_i,\, 1\leqslant i \leqslant h$, let $\psi:P_h \to P_h$ be the linear map defined by $$\psi (g) = t_1^{q_1p^{t+1}+r_1}t_2^{q_2p^{t+1}+r_2}\ldots t_i^{q_ip^{t+1}+r_i}t_{i+1}^{q_{i+1}p^{t+1}}t_{i+2}^{q_{i+2}p^{t+1}}\ldots t_h^{q_hp^{t+1}}g^p$$ for all $g\in P_h$ and $1 \leq r_1, r_2, \ldots, r_i\leq p-1$. Then if $g$ is in the image of $\mathscr A_p$ via $g = \sum_{0\leq k \leq t-2}\mathscr P^k(g_k)$ for some $g_k \in P_h$, then $\psi$ induces a homomorphism
		$$(\mathbb F_p\otimes_{\mathscr A_p} P_h)_n \to (\mathbb F_p\otimes_{\mathscr A_p} P_h)_{pn+2(\sum_{1\leqslant j \leqslant i}(q_jp^{t+1}+r_j)+ \sum_{i+1\leqslant j \leqslant h}q_jp^{t+1})} $$ for all $n \geq 0$.} 
	
	Similar to Theorem 2.6, the statement of this theorem is also not correct. It is well-known that the map $\psi$ induces a homomorphism if and only if $\psi(g)$ is hit for all hit polynomials $g$. Hence, this theorem needs to be restated as follows:
	
	\medskip
	\textit{For positive integers $t \geq 2$, $q_i,\, 1\leqslant i \leqslant h$ and let $\psi:P_h \to P_h$ be the linear map defined by $$\psi (g) = t_1^{q_1p^{t+1}+r_1}t_2^{q_2p^{t+1}+r_2}\ldots t_i^{q_ip^{t+1}+r_i}t_{i+1}^{q_{i+1}p^{t+1}}t_{i+2}^{q_{i+2}p^{t+1}}\ldots t_h^{q_hp^{t+1}}g^p$$ for all $g\in P_h$ and $1 \leq r_1, r_2,\ldots, r_i\leq p-1$. If $g$ is hit, then so is $\psi(g)$. Consequently, the map $\psi$ induces a homomorphism
		$$(\mathbb F_p\otimes_{\mathscr A_p} P_h)_n \to (\mathbb F_p\otimes_{\mathscr A_p} P_h)_{pn+2(\sum_{1\leqslant j \leqslant i}(q_jp^{t+1}+r_j)+ \sum_{i+1\leqslant j \leqslant h}q_jp^{t+1})} $$ for all $n \geq 0$.} 
	
	I predict that the author proved this theorem by using Proposition~3.3 in the paper Tannay-Oner \cite{tan}, but the results in this proposition are only true for some special cases of hit polynomials, hence they do not enough to use for proving this theorem. 
	
	For $h = i = 2,\, q_1 = q_2 = 0,\, r_1 = r_2 = p-1$, this theorem is proved in the paper Crossley \cite[Lemma 2.1, Page 173]{cros} and its proof is long and very complicated. So, I do not believe that the author can prove this theorem.
	
	\item[8.] In the proof of Theorem 2.8 on Page 6, the author stated that 
	``we use an algorithm on MAGMA [3] and find that the invariant
	[$(\mathbb F_p\otimes _{\mathscr A_p}H^*(V;\mathbb F_p))_{n^{(s)}}]^{GL(3,\mathbb F_p)}$ has dimension 1 for every
	$s$.'' 
	
	However, we do not believe that there is a computer program that can check these results for infinitely many values of odd prime numbers $p$. Moreover, this proof is based on unconfirmed results (eg. Theorems 2.4, 2.7) hence, the proof of Theorem 2.8 is refused.
\end{enumerate}

{}

\end{document}